\documentclass[11pt]{article}
\usepackage{setspace}
\usepackage{textcomp}
\usepackage{amsmath, amsthm, amssymb}
\usepackage[left=1.5in,top=1in,right=1in,nohead]{geometry}
\usepackage{hyperref}
\usepackage{latexsym}
\usepackage{amssymb}
\usepackage{mathrsfs}
\usepackage{rawfonts}
\def\beg{\begin}

\usepackage{makeidx}
\makeindex

\def\bct{\begin{center}}
\def\ect{\end{center}}
\def\beg{\begin}
\def\dirac{\rlap{\hspace{.05cm}/}D}
\def\<{\langle}
\def\>{\rangle}
\def\mbb{\mathbb}
\def\mc{\mathcal}
\def\mco{\mathcal O}

\def\mf{\mathfrak}
\def\n{\nabla}

\def\bcp{{\Bbb C\Bbb P}}

\newcounter{exmp}

\title{Deformations of Scalar-Flat Anti-Self-Dual metrics and
Quotients of Enriques Surfaces}

\author{Mustafa Kalafat}
\begin{document}
\maketitle

\begin{abstract}
In this article, we prove that a quotient of a $K3$ surface by a
free ${\mbb Z}_2\oplus{\mbb Z}_2$ action does not admit any metric
of positive scalar curvature. This shows that the scalar flat anti
self-dual metrics (SF-ASD) on this manifold can not be obtained from
a family of metrics for which the scalar curvature changes sign,
contrary to the previously known constructions of this kind of
metrics on manifolds of $b^+=0$
\end{abstract}

%\keywords{Self-Dual Metrics, Spin Structures, Dirac Operator, K\"ahler Manifolds, Algebraic Surfaces}

%\tableofcontents

\section{Introduction}

One of the most interesting features of the space of anti-self-dual
or self-dual(ASD/SD) metrics on a manifold is that the scalar
curvature can change sign on a connected component. That means, one
can possibly join two ASD metrics of scalar curvatures of opposite
signs by a $1$-parameter family of ASD metrics. Whereas, this is not
the case, for example for the space of Einstein metrics. There, each
connected component has a fixed sign for the scalar curvature.

As a consequence, contrary to the Einstein case, most of the
examples of SF-ASD metrics are constructed by first constructing a
family of ASD metrics. Then showing that there are metrics of
positive and negative scalar curvature in the family, and
guaranteeing that there is a scalar-flat member in this family. In
the $b^+=0$ case actually this is the only way known to construct
SF-ASD metrics on a $4$-manifolds. This paper presents an example of
a SF-ASD Riemannian $4$-manifold which is impossible to obtain by
this kind of techniques since it does not have a positive scalar
curvature deformation.

\S\ref{lebrunkim} reviews the known examples of ASD metrics
constructed by a deformation changing the sign of the scalar
curvature, \S\ref{sfasd} introduces an action on the $K3$ surface
and furnish the quotient manifold with a SF-ASD metric,
\S\ref{weitzenbocksection} shows that the smooth manifold defined in
\S\ref{sfasd} does not admit any positive scalar curvature(PSC) or
PSC-ASD metric, finally \S\ref{otherex} includes some related
examples and remarks.

{\bf Acknowledgments.} I want to thank to Claude LeBrun
for his excellent directions. Also I want to thank to %Justin Sawon and
Ioana Suvaina for useful discussions. I would like to express my
gratitude to the organizers of the G\"okova Geometry-Topology
Conference, especially to Turgut \"Onder, Selman Akbulut and Sema
Salur for providing such an academically \'elite and beautiful
scientific atmosphere each year.

\section{Constructions of SF-ASD metrics}\label{lebrunkim}

Here we review some of the constructions for SF-ASD metrics on
$4$-manifolds. We begin with

\beg{thm}[LeBrun\cite{lom}]\label{clsfasd} For all integers $k\geq
6$, the manifold
$$k\overline{\mbb{CP}}_2=\underbrace{\overline{\mbb{CP}}_2\#\cdots\#\overline{\mbb{CP}}_2}_{k-many}$$
admits a $1$-parameter family of real analytic ASD conformal metrics
$[g_t]$ for $t\in [0,1]$ such that $[g_0]$ contains a metric of
$s>0$ on the other hand $[g_1]$ contains a metric of $s<0$
\end{thm}

\beg{cor}[LeBrun\cite{lom}] For all integers $k\geq 6$, the
connected sum ~$k\overline{\mbb{CP}}_2$~ admits scalar-flat
anti-self-dual(SF-ASD) metrics \footnote{Quite recently, LeBrun and
Maskit announced that they have extended this result to the case
$k=5$ with similar techniques, which is the minimal number for these
type of connected sums according to \cite{cltopsd}.}
\end{cor}
\begin{proof}
Let $h_t\in [g_t]$ be a smooth family of metrics representing the
smooth family of conformal classes $[g_t]$ constructed in Lebrun
\cite{lom}.

We know that the smallest eigenvalue $\lambda_t$ of the Yamabe
Laplacian $(\Delta + s/6)$ of the metric $h_t$ exists, and is a
continuous function of $t$. Which measures the sign of the
conformally equivalent constant scalar curvature metric\cite{lp}.

But the theorem(\ref{clsfasd}) tells us that $\lambda_0$ and
$\lambda_1$ has opposite signs. Then there is some $c\in [0,1]$ for
which $\lambda_c=0$. Let $u$ be the eigenfunction corresponding to
the eigenvalue $0$, for the Yamabe Laplacian of $h_c$. Rescale it by
a constant so that it has unit integral.

Since $u$ is a continuous function on the compact manifold, it has a
minimum say at $m$. Choose the normal coordinates around there, so
that $\Delta u(m)=-\sum_{k=1}^4 \partial^2u(m)$.
%by the Remark\ref{normallaplacian}
Second partial derivatives are greater than or equal to zero,
$\Delta u(m)\leq 0$ ~so~ $u(m)=-{6\over s}\Delta u(m)\geq 0$. Assume
$u(m)=0$. Then the maximum of $-u$ is attained and it is nonnegative
with $(-\Delta-s/6)(-u)=0\geq 0$. So the strong maximum
principle %(\ref{strongmax})
is applicable and $-u\equiv 0$. Which is not an eigenfunction. So
$u$ is a positive function. For a conformally equivalent metric, the
change in the scalar curvature is
$$\tilde s=6u^{-3}(\Delta+s/6)u$$
 Thus $g=u^2h_c$ is  a scalar-flat
anti-self-dual metric on $k\overline{\bcp}_2$ for any $k\geq 6$.
\end{proof}

Another construction tells us

\beg{thm}[\cite{kimsfasd}] There exist a continuous family of
self-dual metrics on a connected component of the moduli space of
self-dual metrics on
$$l(S^3\times S^1)\# m \mbb{CP}_2~~\textnormal{for any}~m\geq 1 ~\textnormal{and for
some}~~ l\geq 2$$ which changes the sign of the scalar curvature
\end{thm}

\section{SF-ASD metric on the quotient of Enriques Surface}\label{sfasd}
In this section we are going to describe what we mean by $K3/{\mbb
Z}_2 \oplus {\mbb Z}_2$, and the scalar-flat anti-self-dual(SF-ASD)
metric on it.

Let $A$ and $B$ be real $3\times3$ matrices. For $x,y\in \mbb C^3$,
consider the algebraic variety $V_{2,2,2}\subset \mbb{CP}_5$ given
by the equations
$$\sum_{j}A_i^jx_j^2+B_i^jy_j^2=0~~,~~i=1,2,3$$
more precisely
$$A_1^1x_1^2+A_1^2x_2^2+A_1^3x_3^2+B_1^1x_1^2+B_1^2x_2^2+B_1^3x_3^2=0$$
$$A_2^1x_1^2+A_2^2x_2^2+A_2^3x_3^2+B_2^1x_1^2+B_2^2x_2^2+B_2^3x_3^2=0$$
$$A_3^1x_1^2+A_3^2x_2^2+A_3^3x_3^2+B_3^1x_1^2+B_3^2x_2^2+B_3^3x_3^2=0$$
For generic $A$ and $B$, this is a complete intersection of three
nonsingular quadric hypersurfaces. By the Lefschetz hyperplane
theorem, it is simply connected, and
$$K_{V_2}=K_{\mbb P^5}\otimes [V_2^{\mbb P^5}]=\mco(-6)\otimes\mco(1)^{\otimes2}=\mco(-4) $$
since $[V_2]_h=2[H]_h$ and taking Poincare duals, similarly
$$K_{V_{2,2}}  =K_{V_2}\otimes [V_{2,2}^{\mbb P^5}]=\mco(-4)\otimes\mco(2)=\mco(-2)$$
$$K_{V_{2,2,2}}=K_{V_{2,2}}\otimes[V_{2,2,2}^{\mbb 
P^5}]=\mco(-2)\otimes\mco(2)=\mco$$
finally. So the canonical bundle is trivial. $V$ is a $K3$ Surface.
We define the commuting involutions $\sigma^\pm$ by
$$\sigma^+(x,y)=(x,-y) ~~~~\textnormal{and}~~~~  \sigma^-(x,y)=(\bar x,\bar y)$$
and since we arranged $A$ and $B$ to be real, $\sigma^\pm$ both act
on $V$.\

At a fixed point of $\sigma^+$ on $V$, we have $y_j=-y_j=0$, so
$\sum_j A_i^jx_j^2=0$ implying $\sum_j B_i^jx_j^2=0$, too. So if we
take $A$ and $B$ to be invertible, these conditions are only
satisfied for $x_j=y_j=0$ which does not correspond to a point, so
$\sigma^+$ is free and holomorphic. \ At a fixed point of $\sigma^-$
on $V$, $x_j$'s and $y_j$'s are all real. If $A_1^j , B_1^j > 0$ for
all $j$ then $\sum_{j}A_i^jx_j^2+B_i^jy_j^2=0$ forces $x_j=y_j=0$
making $\sigma^-$ free. \ At a fixed point $\sigma^-\sigma^+$ on
$V$, $x_j=\bar x_j$ and $y_j=-\bar y_j$, so $x_j$'s are real and
$y_j$'s are purely imaginary. Then $y_j^2$ is a negative real
number. So if we choose $A_2^j>0$ and $B_2^j<0$, this forces
$x_j=y_j=0$, again we obtain a free action for $\sigma^-\sigma^+$.
Thus choosing $A$ and $B$ within these circumstances $\sigma^\pm$
generate a free $\mbb Z_2\oplus \mbb Z_2$ action and we define
$K3/{\mbb Z_2\oplus \mbb Z_2}$ to be the quotient of $K3$ by this
free action. We have
$$\chi=\sum_{k=0}^4(-1)^kb_k=2-2b_1+b_2=2+(2b^+-\tau) ~~~\textnormal{hence}~~~ b^+=(\chi+\tau-2)/2$$
so, $b^+(K3/{\mbb Z_2\oplus \mbb Z_2})=(24/4-16/4-2)/2=0$, a special
feature of this manifold.

Next we are going to furnish this quotient manifold with a
Riemannian metric. For that purpose, there is a crucial observation
\cite{hitein} that, for any free involution on $K3$, there exists a
complex structure on $K3$ making this involution holomorphic, so the
quotient is a complex manifold. We begin by stating the

\beg{thm}[Calabi-Yau\cite{calabi,yau,cyghj,joyce}]\index{Calabi-Yau
Theorem}\label{cy} Let $(M,\omega)$ be a compact K\"ahler\\ $n$-manifold.
Given a $(1,1)$-form $\rho$ belonging to the
class $2\pi c_1(M)$ so that it is  closed.\\  Then,\\
there exists a unique K\"ahler metric with form $\omega '$ which is
in the same class as in $\omega$,\\ whose Ricci form is $\rho$
\end{thm}

Intuitively, you can slide the K\"ahler form $\omega$ in its
cohomology class and obtain any desired reasonable Ricci form
$\rho$.

Since $c_1(K3)=0$ in our case, taking $\rho\equiv 0$ gives us a
Ricci-Flat(RF) metric on the $K3$ surface, the {\em Calabi-Yau
metric.} This metric is hyperkahler since, the holonomy group of
K\"ahler manifolds are a subgroup of $U_2$, but Ricci-Flatness
causes a reduction in the holonomy for harmonic forms are parallel
because of the Weitzenb\"ock Formula for spin
manifolds(\ref{weitzenbocks}). Scalar flatness and non triviality of
$b^+$ is to be checked. $b_1(K3)=0$~ implying
~$b^+(K3)=(24-16-2)/2=3$, which is nonzero. Actually $b^+$ is nontrivial for any K\"ahler
surface since the K\"ahler form is harmonic \& self-dual. So we have
the reduction because there are some harmonic parallel forms, the
holonomy group supposed to fix these forms causing a shrinking, the
next possible option is $SU_2$ which is equal to $Sp_1$ in this
dimension, hence the Calabi-Yau metric is hyperk\"ahler.

So we have at least three almost complex structures $I,J,K$,
parallel with respect to the Riemannian connection. By duality we
regard these as three linearly independent self-dual $2$-forms,
parallelizing $\Lambda^+_2$. So any parallel $\Lambda^+_2$ form on
$K3$ defines a complex structure after normalizing. In other words
~$aI+bJ+cK$~ defines a complex structure for the constants
satisfying ~$a^2+b^2+c^2=1$,~ i.e the normalized linear combination.
On the other hand $$b_1(K3/\mbb Z_2)=b_1(K3)=0~~,~~b^+(K3/\mbb
Z_2)=(12-8-2)/2=1$$ Since the pullback of harmonic forms stay
harmonic, the generating harmonic $1$-form on $K3/\mbb Z_2$ is
coming from the universal cover, so is fixed by the $\mbb Z_2$
action. It is also a parallel self-dual form so its normalization is
then a complex structure left fixed by $\mbb Z_2$. So the quotient
is a complex surface with $b_1=0$ and $2c_1=0$ implying that it is
an {\em Enriques Surface.}\index{surface, Enriques}

So we saw that any involution or $\mbb Z_2$-action can be made
holomorphic by choosing the appropriate complex structure on $K3$.
In particular by changing the complex structure, $\sigma^-$ becomes
holomorphic and then both $K3/\mbb Z_2^\pm$ are complex manifolds,
i.e. Enriques Surfaces, for $\mbb Z_2^\pm=\<\sigma^\pm\>$.

Now consider another metric on $K3$ : the restriction of the {\em
Fubini-Study metric} on $\mbb{CP}^5$ obtained from the K\"ahler form
$$\omega_{FS}={i\over 2\pi}\partial\bar\partial \log{|(x_1,x_2,x_3,y_1,y_2,y_3)|^2}$$
We also denote the restriction metric by $g_{FS}$. It is clear that
$\sigma^\pm$ leaves this form invariant, hence they are isometries
of $g_{FS}$. Hence the Fubini-Study metric projects down to the
metrics $g_{FS}^\pm$ on $K3/{\mbb Z_2^\pm}$. Let $h^\pm$ be the
Calabi-Yau metric(\ref{cy}) on $K3/{\mbb Z_2^\pm}$ with K\"ahler
form cohomologous to that of $g_{FS}^\pm$. To remedy the ambiguity
in the negative side, keep in mind that, $\sigma^-$ fixes the metric
and the form on $K3$, though the quotient is not a K\"ahler manifold
initially since it is not a complex manifold, it is locally
K\"ahler. We arrange the complex structure of $K3$ to provide a
complex structure to the form, so the quotient manifold is K\"ahler.
Now we have two K\"ahler metrics on the quotient (for different
complex structures) but we do not know much about their curvatures,
but we want to make the curvature Ricci-Flat, so we use the
Calabi-Yau argument. Since $c_1(K3/\mbb Z_2^\pm)=0$ with real
coefficients, we pass to the Calabi-Yau metric for $\rho=0$.
$\pi^\pm$ denoting the quotient maps, the pullback metrics $\pi^{\pm
*} h^\pm$ are both Ricci-Flat-K\"ahler(RFK) metrics on $K3$ with
K\"ahler forms cohomologous to that of $g_{FS}$. Their Ricci forms
are both $0$. By the uniqueness(\ref{cy}) of the Yau metric we have
$\pi^{+ *} h^+=\pi^{- *} h^-$. Hence this is a Ricci-Flat K\"ahler
metric on $K3$ on which both $\sigma^\pm$ act isometrically. This
metric therefore projects down to a Ricci-Flat metric on our
manifold ~$K3/{\mbb Z_2\oplus \mbb Z_2}$. It is the scalar-flat, and
being locally K\"ahler implies locally scalar-flat anti-self-dual,
hence a SF-ASD metric.

\section{Weitzenb\"ock Formulas}\label{weitzenbocksection}

Now we are going to show that the smooth manifold $K3/{\mbb Z}_2
\oplus {\mbb Z}_2$ does not admit any positive scalar curvature
metric. For that purpose we state the Weitzenb\"ock Formula for the
Dirac Operator on spin manifolds. Before that we introduce some
notation together with some ingredients of the formula.

The Levi-Civita connection is going to be the linear map we denote
by $\nabla : \Gamma(E) \to \Gamma(Hom(TM,E))$ for any vector bundle
$E$ over  a Riemannian Manifold $M$. Then we get the adjoint
 $\nabla^* : \Gamma(Hom(TM,E)) \to \Gamma(E)$ defined implicitly by
 $$\int_M\<\nabla^*S,s\> dV = \int_M\<S,\nabla s\> dV $$  and we
 define the {\em connection Laplacian } \index{Laplacian, connection/rough}
 of a section $s \in \Gamma(E)$ by
 their composition $\nabla^*\nabla s$. Notice that the harmonic
 sections are parallel for this operator. Using the metric, we can
 express its action as :

 \beg{prop}[\cite{petersen}p179] Let $(M,g)$ be an oriented
 Riemannian manifold, $E\to M$ is a vector bundle with an inner
 product and compatible connection. Then
 $$\nabla^*\nabla s= -tr \nabla^2s$$
for all compactly supported sections of $E$
 \end{prop}

\beg{proof} First we need to mention the second covariant
derivatives and then the integral of the divergence.

We set $$\nabla^2K(X,Y)=(\n \n K)(X,Y)=(\n_X \n K)(Y)$$ Then using
the fact that $\n_X$ is a derivation commuting with every
contraction: \cite{kn1}p124 \bct ~ $\n_X \n_Y K=\n_X C(Y \otimes \n
K)=C\n_X(Y \otimes \n K)=C(\n_XY \otimes \n K + Y \otimes \n_X \n
K)=\n_{\n_XY}K + (\n_X \n K)(Y)=\n_{\n_XY}K+\n^2K(X,Y) $ \ect hence
$\n^2K(X,Y)=\n_X \n_Y K - \n_{\n_XY}K$ for any tensor $K$. That is
how the second covariant derivative is defined \index{second
covariant derivative}. Higher covariant derivatives are defined
inductively.

For the divergence, remember that $$(div X) dV_g={\mc L}_XdV_g $$
which is taken as a definition sometimes\cite{kn1}p281. Combining
this with the Cartan's Formula: ${\mc
L}_XdV=di_XdV+i_Xd(dV)=di_XdV$. Then the Stokes' Theorem yields that
$\int_M(div X)dV=\int_M{\mc L}_XdV=\int_Md(i_XdV)=\int_{\partial
M}i_XdV=0$ for a compact manifold without boundary. This is actually
valid even for a noncompact manifold together with a compactly
supported vector field.

Now take an open set on $M$ with an orthonormal basis
$\{E_i\}_{i=1}^n$. Let $s_1$ and $s_2$ be two sections of $E$
compactly supported on the open set. We reduce the left-hand side by
multiplying by $s_2$ as follows: %\bct
\beg{flushleft} $(\n^*\n s_1,s_2)_{L^2}=\int_M\<\n^*\n s_1 ,
s_2\>dV=\int_M\<\n s_1,\n s_2\>dV=\int_Mtr((\n s_1)^*\n s_2)dV$
\end{flushleft} \bct
$=\sum_{i=1}^n \int_M\<(\n s_1)^*\n s_2(E_i),E_i\>dV=\sum\int_M\<(\n
s_1)^*\n_{E_i}s_2,E_i\>dV$ \ect \bct $=\sum \int_M\<\n_{E_i}s_2,\n
s_1(E_i)\>dV=\sum \int_M\<\n_{E_i}s_1 ,\n_{E_i}s_2\>dV$ \ect Define
a vector field $X$ by $g(X,Y)=\<\n_Y s_1,s_2\>$. Divergence of this
vector field is \bct $div X=-d^*(X^\flat)=tr\n X=\sum_{i=1}^n
\<\n_{E_i}X,E_i\>=\sum(E_i\<X,E_i\>-\<X,\n_{E_i}E_i\>)$\ect\bct
$=\sum(E_i\<\n_{E_i}s_1,s_2\>-\<\n_{\n_{E_i}E_i}s_1,s_2\>)$ \ect We
know that its integral is zero, so our expression continues to
evolve as \bct
$\sum\int_M\<\n_{E_i}s_1,\n_{E_i}s_2\>dV-\int_M(divX)dV$ \ect\bct $
=\sum\int_M\<\n_{E_i}s_1,\n_{E_i}s_2\>dV-
\sum\int_M(E_i\<\n_{E_i}s_1,s_2\>-\<\n_{\n_{E_i}E_i}s_1,s_2\>)dV$
\ect\bct
$=\sum\int_M(-\<\n_{E_i}\n_{E_i}s_1,s_2\>+\<\n_{\n_{E_i}E_i}s_1,s_2\>)dV
$ \ect\bct $
=\sum\int_M\<-{\n}^2s_1(E_i,E_i),s_2\>dV=-\int_M\<\sum\<{\n}^2s_1(E_i),E_i\>,s_2\>dV$
\ect\bct $=\int_M\<-tr\n^2s_1,s_2\>dV= (-tr\n^2s_1,s_2)_{L^2} $\ect
So we established that ~$\n^*\n s_1=-tr\n^2 s_1$~ for compactly
supported sections in an open set.

\end{proof}

\beg{thm}[Atiyah-Singer Index
Theorem\cite{lm}p256,\cite{morgansw}p47] \index{index
theorem}\label{indexthm} Let $M$ be a compact spin manifold of
dimension $n=2m$. Then, \\ the index of the Dirac operator is given
by $$
ind({\rlap{\hspace{.05cm}/}D}^+)=\hat{A}(M)=\textbf{\^{A}}(M)[M]$$
More generally, if $E$ is any complex vector bundle over $M$, the
index of \\ ${\rlap{\hspace{.05cm}/}D}^+_E : \Gamma(\mbb
S_\pm\otimes E) \to \Gamma(\Bbb S_\mp\otimes E)$ is given by
$$ind({\rlap{\hspace{.05cm}/}D}^+_E)=\{ch(E)\cdot\textbf{\^{A}}(M) \}[M]$$
For $n=4$,~ $\textbf{\^{A}}(M)=1-p_1/24$ and the first formula
reduces to
$$ ind({\rlap{\hspace{.05cm}/}D}^+)=\hat{A}(M)=\int_M-{p_1(M)\over
24}=-{\tau(M)\over 8}$$ by the Hirzebruch signature Theorem.
%\^A $\mathfrak{A} \hat {\mathcal A}  \hat{\mathsf{A}} \mathrm{A}
%\textsc{A} \hat{\mathtt{A}} \hat{\textbf{\texttt{A}} }
%\mathbf{\hat{\texttt{A}}} $
%~~ \^{\mbb A}
\end{thm}

Let us explain the ingredients beginning by the cohomology class
$\textbf{\^{A}}(M)$. Consider the power series of the following
function\cite{friedrich}p108 :
$${t/2\over \sinh{t/2}}={t\over e^{t/2}-e^{-t/2}}=1+A_2t^2+A_4t^4+\ldots$$
where we compute the coefficients as
$$A_2=-{1\over 24}~~,~~A_4={7\over 10\cdot24\cdot 24}={7\over 5760}$$
Consider the Pontrjagin classes $p_1...p_k$ of $M^{4k}$. Represent
these as the elementary symmetric functions in the squares of the
formal variables $x_1 \cdots x_k$:
$$x_1^2+\cdots +x_k^2=p_1~, ~~\cdots~~ ,~x_1^2x_2^2\cdots x_k^2=p_k$$
Then ~$\prod_{i=1}^k {x_i\over e^{x_i/2}-e^{-x_i/2}} $ ~ is a
symmetric power series in the variables  $x_1^2\cdots x_k^2$, hence
defines a polynomial in the Pontrjagin classes. We call this
polynomial as $\textbf{\^{A}}(M)$
$$\textbf{\^{A}}(M)=\prod_{i=1}^k {x_i/2\over \sinh{x_i/2}}$$
In lower dimensions we have
$$\textbf{\^{A}}(M^4)=1-{1\over 24}p_1~~,~~\textbf{\^{A}}(M^8)=1-{1\over 24}p_1+{7\over 5760}p_1^2-{1\over 1740}p_2$$
If the manifold has dimension $n=4k+2$, again it has $k$ Pontrjagin
classes, and we define the polynomial $\textbf{\^{A}}(M^{4k+2})$ by
the same formulas.

Secondly, we know that ${\rlap{\hspace{.05cm}/}D}^+ : \Gamma(\mbb
S_+) \to \Gamma(\Bbb S_-)$ is an elliptic  % \ref{ellipticity}
operator, so its kernel is finite dimensional and its image is a
closed subspace of finite codimension. The {\em index}\index{index}
of an elliptic operator is defined to be $dim kernel -dim cokernel$.
Actually in our case ${\rlap{\hspace{.05cm}/}D}^+$ and
${\rlap{\hspace{.05cm}/}D}^-$ are formal adjoints of each other:
$(\rlap{\hspace{.05cm}/}D\psi,\eta)_{L^2}=(\psi,\rlap{\hspace{.05cm}/}D\eta)_{L^2}$
for $\psi,\eta$ compactly supported
sections\cite{lm}p114\cite{morgansw}p66. Consequently the index
becomes
$dim ker{\rlap{\hspace{.05cm}/}D}^+ -dim ker{\rlap{\hspace{.05cm}/}D}^-$.\\
This index is computed from the symbol in the following way.
Consider the pullback of ${\mbb S}_\pm$ to the cotangent bundle
$T^*M$. The symbol induces a bundle isomorphism between these
bundles over the complement of the zero section of $T^*M$. In this
way the symbol provides an element in the relative K-theory of
$(T^*M,T^*M-M)$. The Atiyah-Singer Index Theorem computes the index
from this element in the relative K-theory. In the case of the Dirac
operator the index is $\hat{A}(M)$, the so-called {\em A-hat genus}
of $M$.

Now we state our main tool

\beg{thm}[Weitzenb\"ock Formula\cite{petersen}p183,\cite{bes}p55]
 \index{Weitzenb\"ock Formula, spin}\label{weitzenbocks}

On a spin Riemannian manifold, consider the Dirac operator
$\rlap{\hspace{.05cm}/}D : \Gamma(\mbb S_\pm) \to \Gamma(\Bbb
S_\mp)$. The Dirac Laplacian might be expressed in terms of the
connection/rough Laplacian as

$$\rlap{\hspace{.05cm}/}D^2=\nabla^*\nabla+{s\over 4}$$
where $\nabla$ is the Riemannian connection

\end{thm}

Finally we state and prove our main result :

\beg{thm}\label{mainresult}The smooth manifold $K3/{{\mbb Z}_2
\oplus {\mbb Z}_2}$ does not admit any metric of positive scalar
curvature(PSC)
\end{thm}
\beg{proof} If $K3/{{\mbb Z}_2 \oplus {\mbb Z}_2}$ admits a metric
of PSC then $K3$ is also going to admit such a metric because one
pulls back the metric of the quotient, and obtain a locally
identical metric on which the PSC survives.

So we are going to show that the $K3$ surface does not admit any
metric of PSC. First of all the canonical bundle of $K3$ is trivial
so that $c_1(K3)=0=w_2(K3)$ implying that it is a spin manifold.\\
By the Atiyah-Singer Index Theorem \ref{indexthm},
$$ind\rlap{\hspace{.05cm}/}D^+=\textbf{\^{A}}(M)[M]=-{\tau(M)\over 8}=2$$
for the $K3$ Surface. Since it is equal to $dim ker-dim coker$, this
implies that the $dim ker\rlap{\hspace{.05cm}/}D^+\geq 2$. \\ Let
$\varphi \in ker\rlap{\hspace{.05cm}/}D^+$. Then $\dirac
^2\varphi=0$ since $\dirac=\dirac^+\oplus\dirac^-$. So by the spin
Weitzenb\"ock Formula \ref{weitzenbocks}
$$0=\nabla^*\nabla\varphi+{s\over 4}\varphi$$ Taking the inner
product with $\varphi$, so integrating over the manifold yields
$$0=(\nabla^*\nabla\varphi,\varphi)_{L^2}+({s\over
4}\varphi,\varphi)_{L^2}=(\nabla\varphi,\nabla\varphi)_{L^2}+{s\over
4}(\varphi,\varphi)_{L^2}=\int_M(|\nabla\varphi|^2+{s\over
4}|\varphi|^2)dV $$ and $s>0$ implies that
$|\nabla\varphi|=|\varphi|=0$ everywhere, hence $\varphi=0$. So $ker
\dirac^+=0$, which is not the case.\\
Notice that $s\geq0$ and $s(p)>0$ for some point is also enough for
the conclusion because then $\varphi$ would be parallel and zero at
some point implies zero everywhere
\end{proof}

Alternatively, we could use the Weitzenb\"ock Formula for the
Hodge/modern Laplacian to show that there are no PSC
anti-self-dual(ASD) metrics on $K3/{\mbb Z}_2 \oplus {\mbb Z}_2$.
This is a weaker conclusion though sufficient for our purposes

\beg{thm}[Weitzenb\"ock Formula 2\cite{lom}]
 \index{Weitzenb\"ock Formula, Hodge Laplacian}\label{weitzenbockh}

On a Riemannian manifold, the Hodge/modern Laplacian might be
expressed in terms of the connection/rough Laplacian as
$$(d+d^*)^2=\nabla^*\nabla-2W+{s\over 3}$$
where $\nabla$ is the Riemannian connection and $W$ is the Weyl
curvature tensor.
\end{thm}

\beg{thm}The smooth manifold $K3/{{\mbb Z}_2 \oplus {\mbb Z}_2}$
does not admit any anti-self-dual(ASD) metric of positive scalar
curvature(PSC)
\end{thm}

\beg{proof} Again we are going to show this only for $K3$ as in
\ref{mainresult}.  Anti-self-duality reduces our Weitzenbock Formula
\ref{weitzenbockh} to the form
$$(d+d^*)^2=\nabla^*\nabla-2W_-+{s\over 3}$$
because $W=W_-$ or $W_+=0$.

We have already explained that $b_2^+$ of the $K3$ surface is
nonzero. So take a harmonic self-dual $2$-form $\varphi$.
$W_:\Gamma(\lambda^-)\to \Gamma(\lambda^-)$ only acts on
anti-self-dual forms, so it takes $\varphi$ to zero. Applying the
formula
$$0=\nabla^*\nabla\varphi+{s\over 3}\varphi$$
taking the inner product with $\varphi$, so integrating over the
manifold yields similarly
$$0=(\nabla^*\nabla\varphi,\varphi)_{L^2}+({s\over
3}\varphi,\varphi)_{L^2}=(\nabla\varphi,\nabla\varphi)_{L^2}+{s\over
3}(\varphi,\varphi)_{L^2}=\int_M(|\nabla\varphi|^2+{s\over
3}|\varphi|^2)dV $$ and $s>0$ implies that
$|\nabla\varphi|=|\varphi|=0$ everywhere, hence $\varphi=0$. So
$ker( d+d^*)=0$, which is not the case since the space spanned by
the harmonic representatives are already contained. \end{proof}

\section{Other Examples}\label{otherex}

In this section, we will go through some examples. We begin with the
case $b^+=1$.

\beg{thm}[\cite{klp},\cite{rssfk}] For all integers $k\geq 10$, the
connected sum \\ ~$\mbb{CP}_2\#k\overline{\mbb{CP}}_2$~ admits
scalar-flat-K\"ahler(SFK) metrics\footnote{It is a curious fact
that $k=10$ is the minimal number for these type of metrics(SF-ASD)
on $\mbb{CP}_2\#k\overline{\mbb{CP}}_2$, known by \cite{cltopsd}
long before these constructions made. See \cite{lom} for a survey.}
\end{thm}

The case $k\geq 14$ is achieved in \cite{klp}. They start with blow
ups of $\mbb{CP}_1\times \Sigma_2$ the cartesian product of rational
curve and genus-$2$ curve, which already have a SFK metric via the
hyperbolic ansatz of \cite{clexplicit}. After applying an isometric
involution, they got a SFK orbifold, which has isolated
singularities modelled on $\mathbb{C}^2/ \mbb{Z}_2$. Replacing these
singular models with smooth ones, they obtain the desired metric.

For the case $k=10$, Rollin and Singer first construct a related SFK
orbifold with isolated and cyclic singularities of which the algebra
$\mf a_0$ of non-parallel holomorphic vector fields is zero. This is
done by an argument analogous to that of \cite{burnsbart}. The
target manifold is the minimal resolution of this orbifold. To
obtain the target metric, they glue some suitable local models of
SFK metrics to the orbifold. These models are asymptotically locally
Euclidean(ALE) scalar flat K\"ahler metrics constructed in
\cite{calsing}. %They also provide that any blow up of

%Additionally, $K3$ surface has a Hyperk\"ahler metric, could

When a metric is K\"ahler, from the decomposition of the Riemann
Curvature operator, scalar-flatness turns out to be equivalent to
being anti-self-dual. So these metrics are SF-ASD.     %note proof

Since these manifolds have $b^+\neq 0$ Weitzenb\"ock Formulas apply
as in section \S\ref{weitzenbocksection}, so automatically the
scalar curvature can not change sign. These examples show why the
case $b^+=0$ we focussed on, is interesting.

A second type of example is \footnote{Thanks to the referee for pointing out this
example and the remark} \beg{ex} Let $\Sigma_g$ be the genus-$g$
surface with K\"ahler metric of constant curvature $\kappa=-1$, and
$S^2$ be the $2$-sphere with the round $\kappa=+1$ metric. Consider
the product metric on $S^2\times\Sigma_g$ which is K\"ahler with
zero scalar curvature. So it is anti-self-dual. Then we have fixed
point free, orientation reversing, isometric involutions of both
surfaces obtained by antipodal maps. Combination of these
involutions yield an isometric involution on the product and the
metric falls down to a metric on $(S^2\times\Sigma_g)/\mathbb{Z}_2$
which is SF-ASD as these properties are preserved under isometry.
This is an example with all the key properties where the metric is
completely explicit. Note that this  manifold is non-orientable.
\end{ex}

\beg{rmk}The other side of the story discussed here is that we have
$ASD$, conformally flat deformations to negative scalar curvature
metrics. It is obtained by deforming the
$$\rho : \pi_1(M)\longrightarrow SL(2,\mathbb{H})$$
the representation of the fundamental group in $SL(2,\mbb{H})$.

Also, by further investigation, it is possible to get examples which
are doubly covered by e.g. the simply connected examples of
\cite{klp}.
\end{rmk}

\vspace{5mm}
\small \beg{flushleft} \textsc{Dept. of Mathematics, State University of New York, Stonybrook, NY 11790 USA}\\
\textit{E-mail address :} \texttt{\textbf{kalafat@math.sunysb.edu}}
\end{flushleft}

%\printindex


\vspace{5mm}

\newpage






\begin{thebibliography}{70}

\bibitem[Bes]{bes}
{\sc A.~Besse}, {\em {E}instein Manifolds}, Springer-Verlag, 1987.

\bibitem[Burns-Bart]{burnsbart}
{\sc D.~Burns and P.~de Bartolomeis}, {\em Stability of vector
bundles and extremal metrics} Invent. Math. 92, 403--407 (1988)

\bibitem[Ca]{calabi}
{\sc E.~Calabi}, {\em {T}he space of K\"ahler metrics}, In Proceedings
of the ICM, Amsterdam 1954 vol 2 pp.~206--207 North-Holland, Amsterdam
1956

\bibitem[Cal-Sing]{calsing}
{\sc D.~Calderbank and M.~Singer}, {\em Einstein metrics and complex
singularities}, Invent. Math. 156, 405--443 (2004)

\bibitem[Fr]{friedrich}
{\sc T.~Friedrich}, {\em {D}irac operators in Riemannian geometry},
Graduate Studies in Mathematics, 25. American Mathematical Society,
Providence, RI, 2000

\bibitem[GS]{gompf}
{\sc R.~Gompf and A.~Stipsicz}, {\em $4$-Manifolds and Kirby
calculus}, Graduate Studies in Mathematics, 20. American
Mathematical Society, Providence, RI, 1999


\bibitem[GHJ]{cyghj}
{\sc M.~Gross, D.~Huybrechts, D.~Joyce}, {\em Calabi-Yau manifolds
and related geometries}, Lectures from the Summer School held in
Nordfjordeid, June 2001. Universitext. Springer-Verlag, Berlin, 2003


\bibitem[HitEin]{hitein}
{\sc N.~Hitchin}, {\em Compact four-dimensional Einstein manifolds},
 J. Differential Geometry  9  (1974), 435--441


\bibitem[Joyce]{joyce}
{\sc D.~Joyce}, {\em Compact manifolds with special holonomy},
Oxford Mathematical Monographs. Oxford University Press 2000

\bibitem[Kim]{kimsfasd}
{\sc Jongsu ~Kim}, {\em On the scalar curvature of self-dual manifolds},
Math. Ann. 297 (1993), no. 2, 235--251

\bibitem[KimLePon]{klp}
{\sc J.~Kim, C.~LeBrun, M.~Pontecorvo}, {\em Scalar-flat K\"ahler
surfaces of all genera}, J. Reine Angew. Math. 486 (1997), 69--95

\bibitem[KN1]{kn1}
{\sc S.~Kobayashi and K.~Nomizu}, {\em Foundations of differential
geometry 1}, Wiley \& Sons, New York 1963, 1969

\bibitem[LM]{lm}
{\sc H.B.~Lawson and M.L.~Michelsohn}, {\em Spin geometry},
Princeton Mathematical Series, 38. Princeton University Press,
Princeton, NJ, 1989

\bibitem[LeOM]{lom}
{\sc C.~LeBrun},{\em
  Curvature Functionals, Optimal Metrics, and the Differential Topology of
  4-Manifolds}, in Different Faces of Geometry,  Donaldson, Eliashberg,
  and Gromov, editors, Kluwer Academic/Plenum, 2004.

\bibitem[LeExp]{clexplicit}
{\sc C.~LeBrun}, {Explicit self-dual metrics on $C{\rm P}\sb
2\#\cdots\# C{\rm P}\sb 2$} J. Differential Geom. 34 (1991), no. 1,
223--253

\bibitem[LeSD]{cltopsd}
{\sc C.~LeBrun}, {On the Topology of Self-Dual 4-Manifolds}, Proc.
Am. Math. Soc. 98 (1986) 637--640

%\bibitem{lebsing2}
%{\sc C.~LeBrun and M.~Singer}, {\em A {K}ummer-type construction of self-dual
%  {$4$}-manifolds}, Math. Ann., 300 (1994), pp.~165--180.

\bibitem[LP]{lp}
{\sc J.~Lee and T.~Parker}, {\em The {Y}amabe problem}, Bull. Am. Math. Soc.,
  17 (1987), pp.~37--91.


\bibitem[Mci]{mcinnes}
{\sc Brett ~McInnes}, {\em Methods of holonomy theory for Ricci-flat
Riemannian manifolds}, J. Math. Phys. 32 (1991), no. 4, 888--896


\bibitem[MoSW]{morgansw}
{\sc J.W.~Morgan}, {\em The Seiberg-Witten equations and
applications to the topology of smooth four-manifolds}, Mathematical
Notes, 44. Princeton University Press, Princeton, NJ, 1996


\bibitem[Pet]{petersen}
{\sc P.~Petersen}, {\em Riemannian geometry},
  Graduate Texts in Mathematics, 171. Springer-Verlag, New York, 1998


\bibitem[Poon86]{poon1}
{\sc Y.S.~Poon}, {\em Compact self-dual manifolds with positive scalar curvature},
J. Differential Geom. 24 (1986), no. 1, 97--132

\bibitem[RS-SFK]{rssfk}
{\sc  Y.~Rollin, M.~Singer}, {\em Non-minimal scalar-flat K\"ahler
surfaces and parabolic stability}, Invent. Math. 162 (2005), no. 2,
235--270

\bibitem[Yau]{yau}
{\sc S.T.~Yau}, {\em On the Ricci curvature of a compact K\"ahler
manifold and the complex Monge-Amp\'ere equation I}, Comm. Pure Appl.
Math. 31 (1978), no. 3, 339--411


\end{thebibliography}
  \end{document}